\documentclass{amsart}

\usepackage{fancyhdr}
\usepackage{lastpage}
\usepackage{stmaryrd,yhmath}

\newcommand{\bdis}{\begin{displaymath}}
\newcommand{\edis}{\end{displaymath}}
\newcommand{\be}{\begin{equation}}
\newcommand{\ee}{\end{equation}}
\newcommand{\mbb}{\mathbb}
\newcommand{\mcal}{\mathcal}

\newcommand{\vp}{\varphi}

\newcommand{\vth}{\vartheta}

\newcommand{\mT}{\mathring{T}}

\newcommand{\zf}{\zeta\left(\frac{1}{2}+it\right)}

\newtheorem{lemma}[]{Lemma}

\theoremstyle{definition}

\theoremstyle{remark}
\newtheorem{remark}[]{Remark}

\newtheorem*{mydef1}{{\bf Theorem}}

\newtheorem*{mydef4}{{\bf Corollary}}

\newtheorem*{mydef81}{{\bf Property 1}}

\newtheorem*{mydef82}{{\bf Property 2}}

\numberwithin{equation}{section}

\begin{document}

\title{Jacob's ladders, reverse iterations and new infinite set of $L_2$-orthogonal systems generated by the Riemann
$\zf$-function}

\author{Jan Moser}

\address{Department of Mathematical Analysis and Numerical Mathematics, Comenius University, Mlynska Dolina M105, 842 48 Bratislava, SLOVAKIA}

\email{jan.mozer@fmph.uniba.sk}

\keywords{Riemann zeta-function}

\begin{abstract}
It is proved in this paper that continuum set of $L_2$-orthogonal systems generated by the Riemann zeta-function on the critical line corresponds to
every fixed $L_2$-orthogonal system on a fixed segment. This theorem serves as a resource for new set of integrals not accessible by the current
methods in the theory of the Riemann zeta-function. \\

\noindent
Dedicated to the 100th anniversary of G.H. Hardy's fundamental theorem: the function $\zf$ has an infinite set of zeros, \cite{1}.
\end{abstract}
\maketitle

\section{Introduction}

\subsection{}

In this paper we obtain new properties of the signal
\be \label{1.1}
\begin{split}
 & Z(t)=e^{i\vth(t)}\zf, \\
 & \vth(t)=-\frac t2\ln \pi +\text{Im}\ln\Gamma\left(\frac 14+i\frac t2\right),
\end{split}
\ee
which is generated by the Riemann zeta-function. In connection with (\ref{1.1}) we have introduced (see \cite{5}, (9.1), (9.2)) the formula
\be \label{1.2}
\tilde{Z}^2(t)=\frac{{\rm d}\vp_1(t)}{{\rm d}t},
\ee
where
\be \label{1.3}
\begin{split}
& \tilde{Z}^2(t)=\frac{Z^2(t)}{2\Phi'_\vp[\vp(t)]}=\frac{\left|\zf\right|^2}{\omega(t)}, \\
& \omega(t)=\left\{ 1+\mcal{O}\left(\frac{\ln\ln t}{\ln t}\right)\right\}\ln t.
\end{split}
\ee
The function $\vp_1(t)$ that we call Jacob's ladder (see our paper \cite{2}) according to the Jacob's dream in Chumash, Bereishis, 28:12, has
the following properties:
\begin{itemize}
\item[(a)]
\bdis
\vp_1(t)=\frac 12\vp(t),
\edis
\item[(b)] function $\vp(t)$ is solution of the non-linear integral equation (see \cite{2}, \cite{5})
\bdis
\int_0^{\mu[x(T)]}Z^2(t)e^{-\frac{2}{x(T)}t}{\rm d}t=\int_0^T Z^2(t){\rm d}t,
\edis
where each admissible function $\mu(y)$ generates a solution
\bdis
y=\vp_\mu(T)=\vp(T);\ \mu(y)\geq 7y\ln y.
\edis
\end{itemize}

\begin{remark}
The main reason to introduce Jacob's ladders in \cite{2} lies in the following: the Hardy-Littlewood integral (1918)
\bdis
\int_0^T \left|\zf\right|^2{\rm d}t
\edis
has -- in addition to the Hardy-Littleweood (and other similar) expression possessing unbounded errors at $T\to\infty$ -- the following infinite set
of almost exact expressions
\be \label{1.4}
\begin{split}
& \int_0^T \left|\zf\right|^2{\rm d}t=\vp_1(T)\ln\vp_1(T)+(c-\ln2\pi)\vp_1(T)+ \\
& + c_0+\mcal{O}\left(\frac{\ln T}{T}\right),\ T\to\infty,
\end{split}
\ee
where $c$ is the Euler's constant, and $c_0$ is the constant from the Titchmarsh-Kober-Atkinson formula (see \cite{7}, p. 141).
\end{remark}

\begin{remark}
Simultaneously with (\ref{1.4}) we have proved that the following transcendental equation
\bdis
\int_0^T \left|\zf\right|^2{\rm d}t=V(T)\ln V(T)+(c-\ln 2\pi)V(T)+c_0
\edis
has an infinite set of asymptotic solutions
\bdis
V(T)=\vp_1(T),\ T\to\infty.
\edis
\end{remark}

\begin{remark}
The Jacob's ladder $\vp_1(T)$ can be interpreted by our formula (see \cite{2}, (6.2))
\be \label{1.5}
T-\vp_1(T)\sim (1-c)\pi(T); \ \pi(T)\sim \frac{T}{\ln T},\ T\to\infty,
\ee
where $\pi(T)$ is the prime-counting function, as an asymptotic complement function to the function
\bdis
(1-c)\pi(T)
\edis
in the sense
\bdis
\vp_1(T)+(1-c)\pi(T)\sim T,\ T\to\infty.
\edis
\end{remark}

\subsection{}

In the paper \cite{3} we have proved that the following continuum set $S(T,2l)$ of the systems
\be \label{1.6}
\begin{split}
& \left\{ |\tilde{Z}(t)|,|\tilde{Z}(t)|\cos\left[\frac{\pi}{l}(\vp_1(t)-T)\right],|\tilde{Z}(t)|\sin\left[\frac{\pi}{l}(\vp_1(t)-T)\right], \dots \right. \\
& \left. |\tilde{Z}(t)|\cos\left[\frac{\pi}{l}n(\vp_1(t)-T)\right],|\tilde{Z}(t)|\sin\left[\frac{\pi}{l}n(\vp_1(t)-T)\right],\dots\right\}, \\
& t\in [\overset{1}{T},\overset{1}{\wideparen{T+2l}}]; \\
& \vp_1\left\{[\overset{1}{T},\overset{1}{\wideparen{T+2l}}]\right\}=[T,T+2l]
\end{split}
\ee
is the set of orthogonal system on the segment
\bdis
[\overset{1}{T},\overset{1}{\wideparen{T+2l}}]
\edis
for all
\bdis
T\geq T_0[\vp_1],\quad 2l\in \left.\left( 0, \frac{T}{\ln T}\right]\right. .
\edis
Next, in the paper \cite{4} we have constructed corresponding continuum set of orthogonal systems generated by Jacobi's polynomials. \\

In this paper we give essential generalization of above mentioned. Namely, to every fixed $L_2$-orthogonal system
\bdis
\{ f_n(t)\}_{n=1}^\infty,\quad t\in [0,2l]
\edis
we assign continuum set of $L_2$-orthogonal systems
\bdis
\begin{split}
& \{ F_n(t;T,k,l)\}_{n=1}^\infty,\ t\in [\overset{k}{T},\overset{k}{\wideparen{T+2l}}],\ T\to\infty,\ k=1,\dots,k_0, \\
& l=o\left(\frac{T}{\ln T}\right); \\
& \vp_1\left\{[\overset{k}{T},\overset{k}{\wideparen{T+2l}}]\right\}=[\overset{k-1}{T},\overset{k-1}{\wideparen{T+2l}}],
\end{split}
\edis
where $k_0\in \mbb{N}$ is an arbitrary fixed number.

\section{Result}

\subsection{}

Let us remind that (see \cite{6})
\bdis
\vp_1^r(t):\ \vp_1^0(t)=t,\ \vp_1^1(t)=\vp_1(t),\ \vp_1^2(t)=\vp_1(\vp_1(t)), \dots
\edis
The following Theorem holds true.

\begin{mydef1}
For every fixed $L_2$-orthogonal system
\be \label{2.1}
\{ f_n(t)\}_{n=1}^\infty,\quad t\in [0,2l],\ l=o\left(\frac{T}{\ln T}\right),\ T\to\infty
\ee
there is continuum set of $L_2$-orthogonal systems
\be \label{2.2}
\begin{split}
& \{ F_n(t;T,k,l)\}_{n=1}^\infty=\\
& = \left\{ f_n(\vp_1^k(t)-T)\prod_{r=0}^{k-1}\left|\tilde{Z}[\vp_1^r(t)]\right|\right\}_{n=1}^\infty, \
t\in [\overset{k}{T},\overset{k}{\wideparen{T+2l}}],
\end{split}
\ee
where
\be \label{2.3}
\begin{split}
& \vp_1\left\{[\overset{k}{T},\overset{k}{\wideparen{T+2l}}]\right\}=[\overset{k-1}{T},\overset{k-1}{\wideparen{T+2l}}],\ k=1,\dots,k_0, \\
& [\mT,\widering{T+2l}]=[T,T+2l],\ T\to\infty
\end{split}
\ee
and $k_0\in\mbb{N}$ is arbitrary number, i. e.  the following formula is valid
\be \label{2.4}
\begin{split}
& \int_{\overset{k}{T}}^{\overset{k}{\wideparen{T+2l}}}f_m(\vp_1^k(t)-T)f_n(\vp_1^k(t)-T)\prod_{r=0}^{k-1}
\tilde{Z}^2[\vp_1^r(t)]{\rm d}t=\\
& = \left\{\begin{array}{rcl} 0 & , & m\not= n, \\ A_n & , & m=n,  \end{array}\right.\quad
A_n=\int_0^{2l} f_n^2(t){\rm d}t.
\end{split}
\ee
Next, we have the following properties
\bdis
l=o\left(\frac{T}{\ln T}\right) \ \Rightarrow
\edis
\be \label{2.5}
|[\overset{k}{T},\overset{k}{\wideparen{T+2l}}]|=\overset{k}{\wideparen{T+2l}}-\overset{k}{T}=o\left(\frac{T}{\ln T}\right),
\ee
\be \label{2.6}
|[\overset{k-1}{\wideparen{T+2l}},\overset{k}{T}]|=\overset{k}{T}-\overset{k-1}{\wideparen{T+2l}}\sim
(1-c)\pi(T);\ \pi(T)\sim \frac{T}{\ln T},
\ee
\be \label{2.7}
[T,T+2l]\prec [\overset{1}{T},\overset{1}{\wideparen{T+2l}}]\prec \dots \prec
[\overset{k}{T},\overset{k}{\wideparen{T+2l}}]\prec \dots ,
\ee
where $\pi(T)$ stands for the prime-counting function.
\end{mydef1}

\subsection{}

\begin{remark}
We obtain from (\ref{2.2}) by (\ref{1.3}) that
\be \label{2.8}
F_n(t;T,k,l)=f_n(\vp_1^k(t)-T)\prod_{r=0}^{k-1}
\frac{\left|\zeta\left(\frac 12+i\vp_1^r(t)\right)\right|}{\sqrt{\omega[\vp_1^r(t)]}},
\ee
i. e. our formula (\ref{2.8}) shows direct connection between the Riemann function
\bdis
\zf
\edis
and an arbitrary $L_2$-orthogonal system
\bdis
\{ f_n(t)\}_{n=1}^\infty,\ t\in [0,2l].
\edis
\end{remark}

\begin{remark}
Asymptotic behavior of the disconnected set (see (\ref{2.6}), (\ref{2.7}))
\be \label{2.9}
\Delta(T,k,l)=\bigcup_{r=0}^k [\overset{r}{T},\overset{r}{\wideparen{T+2l}}]
\ee
is as follows: if $T\to\infty$, then the components of the set (\ref{2.9}) recedes unboundedly each from other
and all together are receding to infinity. Hence, if $T\to\infty$ the set (\ref{2.9}) behaves as one
dimensional Friedmann-Hubble expanding universe.
\end{remark}

\subsection{}

Since (see (\ref{2.3}))
\bdis
\begin{split}
 & t\in [\overset{k}{T},\overset{k}{\wideparen{T+2l}}] \ \Rightarrow \\
 & \vp_1(t)\in [\vp_1(\overset{k}{T}),\vp_1(\overset{k}{\wideparen{T+2l}})]=
 [\overset{k-1}{T},\overset{k-1}{\wideparen{T+2l}}] \ \Rightarrow \\
 & \vp_1^2(t)\in [\vp_1(\overset{k-1}{T}),\vp_1(\overset{k-1}{\wideparen{T+2l}})]=
 [\overset{k-2}{T},\overset{k-2}{\wideparen{T+2l}}] \ \Rightarrow \\
 & \vdots
\end{split}
\edis
we point-out the following

\begin{mydef81}
If
\bdis
t\in [\overset{k}{T},\overset{k}{\wideparen{T+2l}}],\ k=1,\dots,k_0
\edis
then
\be \label{2.10}
\vp_1^r(t)\in [\overset{k-r}{T},\overset{k-r}{\wideparen{T+2l}}],\ r=0,1,\dots,k
\ee
holds true for the arguments of the functions (see (\ref{2.2}), (\ref{2.8}))
\bdis
f_n(\vp_1^k(t)-T), |\tilde{Z}[\vp_1^r(t)]|, \omega[\vp_1^r(t)],
\left|\zeta\left(\frac 12+i\vp_1^r(t)\right)\right|.
\edis
\end{mydef81}

\section{Examples}

\subsection{}

For the classical Fourier orthogonal system
\be \label{3.1}
\begin{split}
 & \left\{ 1,\cos\frac{\pi t}{l},\sin\frac{\pi t}{l},\dots,\cos\frac{\pi n t}{l},\cos\frac{\pi n t}{l},
 \dots\right\} \\
 & t\in [0,2l],\ l=o\left(\frac{T}{\ln T}\right)
\end{split}
\ee
we have as corresponding (see (\ref{2.2}), (\ref{2.8})) continuous set of orthogonal systems the following
\be \label{3.2}
\begin{split}
 & \left\{ \prod_{r=0}^{k-1}\frac{\left|\zeta\left(\frac 12+i\vp_1^r(t)\right)\right|}{\sqrt{\omega[\vp_1^r(t)]}},
 \dots ,\right. \\
 & \left. \left(\prod_{r=0}^{k-1}\frac{\left|\zeta\left(\frac 12+i\vp_1^r(t)\right)\right|}{\sqrt{\omega[\vp_1^r(t)]}}\right)\cos\left(\frac{\pi}{l}n(\vp_1^k(t)-T)\right), \right. \\
 & \left. \left(\prod_{r=0}^{k-1}\frac{\left|\zeta\left(\frac 12+i\vp_1^r(t)\right)\right|}{\sqrt{\omega[\vp_1^r(t)]}}\right)\sin\left(\frac{\pi}{l}n(\vp_1^k(t)-T)\right), \dots \right\} , \\
 & t\in [\overset{k}{T},\overset{k}{\wideparen{T+2l}}],\ k=1,\dots,k_0, \ T\to\infty,
\end{split}
\ee
and, for example,
\bdis
k_0=S=10^{10^{10^{34}}},\ S^S,\dots
\edis
where $S$ is the Skeewes constant.

\subsection{}

For the system of Jacobi's functions
\be \label{3.3}
\sqrt{(1-t)^\alpha(1+t)^\beta}P^{(\alpha,\beta)}_{n}(t),\ t\in [-1,1],\ n=0,1,2,\dots; \
\alpha,\beta>-1
\ee
generated by the Jacobi's polynomials $P^{(\alpha,\beta)}_{n}$ we have that
\be \label{3.4}
\begin{split}
 & \int_{-1}^1 (1-t)^\alpha(1+t)^\beta P^{(\alpha,\beta)}_{n}(t)P^{(\alpha,\beta)}_{m}{\rm d}t=0,\ m\not=n, \\
 & \int_{-1}^1 (1-t)^\alpha(1+t)^\beta \left[P^{(\alpha,\beta)}_{n}(t)\right]^2{\rm d}t= \\
 & = \frac{2^{\alpha+\beta+1}}{2n+\alpha+\beta+1}
 \frac{\Gamma(n+\alpha+1)\Gamma(n+\beta+1)}{n!\Gamma(n+\alpha+\beta+1)}=A_n(\alpha,\beta).
\end{split}
\ee
Next, the substition
\bdis
x=t-T-1
\edis
in (\ref{3.3}) yields (see (\ref{3.4})) the formulae
\bdis
\int_T^{T+2}(2+T-t)^\alpha(t-T)^\beta P^{(\alpha,\beta)}_{m}(t-T-1)P^{(\alpha,\beta)}_{n}(t-T-1){\rm d}t=0,\
m\not=n, \dots
\edis
Consequently, the following continuum set (for each fixed pair $\alpha,\beta>-1$) of orthogonal systems
\bdis
\begin{split}
 & \left\{ P^{(\alpha,\beta)}_{n}(t-T-1)\sqrt{(T+2-\vp_1^k(t))^\alpha(\vp_1^k(t)-T)^\beta} \times \prod_{r=0}^{k-1}
 \frac{\left|\zeta\left(\frac 12+i\vp_1^r(t)\right)\right|}{\sqrt{\omega[\vp_1^r(t)]}}\right\}_{n=0}^\infty, \\
 & t\in [\overset{k}{T},\overset{k}{\wideparen{T+2}}],\ T\to\infty,\ k=1,\dots,k_0 .
\end{split}
\edis
corresponds to the Jacobi's orthogonal system (\ref{3.3}) (see (\ref{2.2}), (\ref{2.8})).

\subsection{}

For the system of Bessel's functions
\be \label{3.5}
\left\{ \sqrt{t}J_n\left(\frac{\mu_m^{(n)}}{2l}t\right)\right\}_{m=1}^\infty,\ t\in [0,2l]
\ee
generated by Bessel's function $J_n(t)$ we have that
\bdis
\begin{split}
 & \int_0^{2l} t J_n\left(\frac{\mu_{m_1}^{(n)}}{2l}t\right)J_n\left(\frac{\mu_{m_2}^{(n)}}{2l}t\right){\rm d}t=0,\ m_1\not=m_2, \\
 & \int_0^{2l} t \left[J_n\left(\frac{\mu_{m}^{(n)}}{2l}t\right)\right]^2{\rm d}t=
 2l^2\left[J_n'(\mu_m^{(n)})\right]^2,
\end{split}
\edis
where
\bdis
\{ \mu_m^{(n)}\}_{m=1}^\infty
\edis
is the sequence of the roots of equation
\bdis
J_n(\mu)=0.
\edis
Consequently, the following continuum set (for each fixed $n$) of orthogonal systems
\bdis
\begin{split}
 &
 \left\{ J_n\left(\frac{\mu_m^{(n)}}{2l}(\vp_1^k(t)-T)\right)\sqrt{\vp_1^k(t)-T}
 \prod_{r=0}^{k-1}
 \frac{\left|\zeta\left(\frac 12+i\vp_1^r(t)\right)\right|}{\sqrt{\omega[\vp_1^r(t)]}}\right\}_{m=1}^\infty , \\
 & t\in [\overset{k}{T},\overset{k}{\wideparen{T+2l}}],\ T\to\infty,\ k=1,\dots,k_0 .
\end{split}
\edis
corresponds to the Bessel orthogonal system (\ref{3.5}) (see (\ref{2.2}), (\ref{2.8})).

\section{Formula (\ref{2.4}) as a resource of new integrals containing multiples of $|\zeta|^2$}

We consider the formula (see (\ref{2.4}), (\ref{2.8}))
\be \label{4.1}
\begin{split}
 & \int_{\overset{k}{T}}^{\overset{k}{\wideparen{T+2l}}}
 f_n^2(\vp_1^k(t)-T)\prod_{r=0}^{k-1}
 \frac{\left|\zeta\left(\frac 12+i\vp_1^r(t)\right)\right|^2}{\omega[\vp_1^r(t)]}{\rm d}t=A_n, \\
 & A_n=\int_0^{2l} f_n^2(t){\rm d}t,\ n=1,2,\dots
\end{split}
\ee

\subsection{}

Let (see (\ref{2.9})
\bdis
t\in \Delta^0(T,k,l)=\bigcup_{r=0}^k (\overset{r}{T},\overset{r}{\wideparen{T+2l}}),\
k=1,\dots,k_0.
\edis
Of course
\bdis
\Delta^0(T,k,l)\subset [T,\overset{k}{\wideparen{T+2l}}],
\edis
and (see (\ref{2.5}) -- (\ref{2.7}))
\be \label{4.2}
\begin{split}
 & |[T,\overset{k}{\wideparen{T+2l}}]|=\sum_{r=0}^k |[\overset{r}{T},\overset{r}{\wideparen{T+2l}}]|+
 \sum_{r=1}^k|[\overset{r-1}{\wideparen{T+2l}},\overset{r}{T}]|= \\
 & =(k+1)o\left(\frac{T}{\ln T}\right)+k\mcal{O}\left(\frac{T}{\ln T}\right)= \\
 & = \mcal{O}\left(\frac{T}{\ln T}\right),\ k=1,\dots,k_0.
\end{split}
\ee
Thus, we have the following: if
\bdis
t\in [T,\overset{k}{\wideparen{T+2l}}],
\edis
then (see (\ref{4.2}))
\bdis
\ln t=\ln(t-T+T)=\ln T+\ln\left(1+\frac{t-T}{T}\right)=\ln T+ \mcal{O}\left(\frac{1}{\ln T}\right),
\edis
i. e.
\be \label{4.3}
\ln t\sim\ln T,\ \forall\- t\in (T,\overset{k}{\wideparen{T+2l}}),\ k=1,\dots,k_0.
\ee

\subsection{}

It is sufficient to use, for example, the formula (\ref{4.1}) in the case (see (\ref{3.1}))
\bdis
f(t)=1,\ \Rightarrow \ A_1=2l ,
\edis
i. e.
\be \label{4.4}
\int_{\overset{k}{T}}^{\overset{k}{\wideparen{T+2l}}}
\prod_{r=0}^{k-1}
 \frac{\left|\zeta\left(\frac 12+i\vp_1^r(t)\right)\right|^2}{\omega[\vp_1^r(t)]}{\rm d}t=2l.
\ee
Next, we obtain from (\ref{4.4}) (see (\ref{1.3}), (\ref{4.3})) by the mean-value theorem that
\be \label{4.5}
\int_{\overset{k}{T}}^{\overset{k}{\wideparen{T+2l}}}
\prod_{r=0}^{k-1}\left|\zeta\left(\frac 12+i\vp_1^r(t)\right)\right|^2 {\rm d}t\sim 2l\ln^kT,\ T\to\infty.
\ee
Consequently, we obtain from (\ref{4.5}) in the case
\bdis
2l=\frac{\Omega}{\ln^kT}=o\left(\frac{T}{\ln T}\right),\ \Omega>0
\edis
the following

\begin{mydef4}
\be \label{4.6}
\int_{\overset{k}{T}}^{\overset{k}{\wideparen{T+\Omega\ln^{-k}T}}}
\prod_{r=0}^{k-1}\left|\zeta\left(\frac 12+i\vp_1^r(t)\right)\right|^2 {\rm d}t\sim \Omega,\ T\to\infty ,
\ee
where
\bdis
0<\Omega=o(T\ln^{k-1}T),\ k=1,\dots,k_0.
\edis
\end{mydef4}

\begin{remark}
Let us notice explicitly that nor the first two formulae (see (\ref{4.6}), $k=1,2;\ \Omega=1$)
\be \label{4.7}
\begin{split}
 & \int_{\overset{1}{T}}^{\overset{1}{\wideparen{T+\ln^{-1}T}}}\left|\zf\right|^2{\rm d}t\sim 1, \\
 & \int_{\overset{1}{T}}^{\overset{1}{\wideparen{T+\ln^{-1}T}}}\left|\zf\right|^2
 \left|\zeta\left(\frac 12+i\vp_1(t)\right)\right|^2{\rm d}t\sim 1,\ T\to\infty
\end{split}
\ee
are not accessible by the current methods in the theory of the Riemann zeta-function.
\end{remark}

\begin{remark}
The first formula in (\ref{4.7}) gives us the answer to the question about a form of segments for which the following
\bdis
\begin{split}
 & [a(T),b(T)] \ \rightarrow \\
 & \int_{a(T)}^{b(T)} \left|\zf\right|^2{\rm d}t\sim 1,\ T\to\infty
\end{split}
\edis
holds true. Namely, corresponding segments are as follows
\bdis
[a(T),b(T)]=[\overset{1}{T},\overset{1}{\wideparen{T+\frac{1}{\ln T}}}]\leftarrow [T,T+\frac{1}{\ln T}].
\edis
\end{remark}

\section{First lemmas}

\subsection{}

The sequence
\be \label{5.1}
\{ \overset{k}{T}\}_{k=1}^\infty
\ee
is defined by the formula (comp. (\ref{2.3}))
\be \label{5.2}
\vp_1(\overset{k}{T})=\overset{k-1}{T},\ k=1,\dots,k_0,\ \overset{0}{T}=T
\ee
for every $T\geq T_0[\vp_1]$, where $k_0\in\mbb{N}$ is an arbitrary fixed number. Since the function
\bdis
\vp_1(t),\ t\to\infty
\edis
increases to $\infty$, then we have from (\ref{5.2}) that
\bdis
\left\{ \overset{k-1}{T}\to\infty\right\} \ \Leftrightarrow \ \left\{ \overset{k}{T}\to\infty\right\},
\edis
i. e.
\be \label{5.3}
\left\{ T\to\infty \right\} \ \Leftrightarrow \ \left\{ \overset{k}{T}\to\infty\right\}.
\ee
Next, we have (see (\ref{1.5}), (\ref{5.2}), (\ref{5.3})) that
\be \label{5.4}
\overset{k}{T}-\overset{k-1}{T}\sim (1-c)\frac{\overset{k}{T}}{\ln \overset{k}{T}} \ \Rightarrow \
1-\frac{\overset{k-1}{T}}{\overset{k}{T}}\sim \frac{1-c}{\ln T},
\ee
i. e.
\bdis
\overset{k-1}{T}\sim \overset{k}{T},
\edis
and, consequently,
\be \label{5.5}
\overset{k}{T}=\{ 1+o(1)\}T,\ T\to\infty,\ k=1,\dots,k_0.
\ee
Since
\bdis
\frac{\overset{k}{T}}{\ln \overset{k}{T}}\sim \frac{T}{\ln T},
\edis
(see (\ref{5.5})) then we have (see {\ref{5.4})) for the sequence (\ref{5.1}) that
\be \label{5.6}
\overset{k}{T}-\overset{k-1}{T}\sim (1-c)\frac{T}{\ln T},\ T\to\infty,\ k=1,\dots,k_0.
\ee
Consequently, we have
\be \label{5.7}
T< \overset{1}{T}<\dots <\overset{k_0}{T}
\ee
and, of course, (see (\ref{2.1}))
\be \label{5.8}
T+H<\overset{1}{\wideparen{T+H}}<\dots <\overset{k_0}{\wideparen{T+H}},\
0<H=o\left(\frac{T}{\ln T}\right).
\ee

\subsection{}

The following lemma holds true.

\begin{lemma}
\be \label{5.9}
\begin{split}
 & H=o\left(\frac{T}{\ln T}\right) \ \Rightarrow \\
 & |[\overset{k}{T},\overset{k}{\wideparen{T+H}}]|=\overset{k}{\wideparen{T+H}}-\overset{k}{T}=
 o\left(\frac{T}{\ln T}\right),\ T\to\infty, \ k=1,\dots,k_0,
\end{split}
\ee
i. e. (\ref{2.5}) holds true.
\end{lemma}

\begin{proof}
 First of all, it follows from (\ref{5.6}) that
\bdis
\overset{k}{T}-T\sim (1-c)k\frac{T}{\ln T},
\edis
i. e.
\be \label{5.10}
\overset{k}{T}-T=\{ 1+o_1(1)\}(1-c)k\frac{T}{\ln T}
\ee
and, simultaneously (see (\ref{5.8}))
\be \label{5.11}
\overset{k}{\wideparen{T+H}}-(T+H)=\{ 1+o_2(1)\}(1-c)k\frac{T}{\ln T} .
\ee
Then we have (see (\ref{5.9}) -- (\ref{5.11})) that
\bdis
\begin{split}
 & 0<\overset{k}{\wideparen{T+H}}-\overset{k}{T}=H+[o_2(1)-o_1(1)](1-c)k\frac{T}{\ln T}= \\
 & = H+[o_4(1)-o_3(1)]\frac{T}{\ln T} = \\
 & = o\left(\frac{T}{\ln T}\right)+o(1)\frac{T}{\ln T}=\\
 & = o\left(\frac{T}{\ln T}\right),\ T\to\infty.
\end{split}
\edis
\end{proof}

\subsection{}

Next, the following lemma holds true
\begin{lemma}
\be \label{5.12}
\begin{split}
 & H=o\left(\frac{T}{\ln T}\right) \ \Rightarrow \\
 & \overset{k}{T}-\overset{k-1}{\wideparen{T+H}}\sim (1-c)\frac{T}{\ln T},\ T\to\infty,\
 k=1,\dots,k_0,
\end{split}
\ee
i. e. (\ref{2.6}) holds true.
\end{lemma}

\begin{proof}
We have from (\ref{5.6}) by (\ref{5.8}), (\ref{5.9}) that
\bdis
\begin{split}
 & \overset{k}{T}-\overset{k-1}{\wideparen{T+H}}+\overset{k-1}{\wideparen{T+H}}-\overset{k-1}{\wideparen{T}}
 \sim (1-c)\frac{T}{\ln T}, \\
 & \overset{k}{T}-\overset{k-1}{\wideparen{T+H}}\sim (1-c)\frac{T}{\ln T}-
 (\overset{k-1}{\wideparen{T+H}}-\overset{k-1}{\wideparen{T}})\sim \\
 & \sim (1-c)\frac{T}{\ln T}+o\left(\frac{T}{\ln T}\right)\sim \\
 & \sim (1-c)\frac{T}{\ln T},\ T\to\infty,\ k=1,\dots,k_0.
\end{split}
\edis
\end{proof}

\begin{remark}
We have (see (\ref{5.12})) that
\be \label{5.13}
[T,T+H]\prec [\overset{1}{T},\overset{1}{\wideparen{T+H}}]\prec \dots \prec
[\overset{k_0}{T},\overset{k_0}{\wideparen{T+H}}],
\ee
i. e. (\ref{2.7}) holds true.
\end{remark}

\section{Reverse iterations}

\subsection{}

First of all, we have (see (\ref{2.3}), (\ref{5.2})) that
\be \label{6.1}
\vp_1(\overset{k}{T})=\overset{k-1}{T} \ \Rightarrow \dots \Rightarrow
\vp_1^k(\overset{k}{T})=T,\ k=1,\dots,k_0.
\ee
Since
\be \label{6.2}
\vp_1(\overset{1}{T})=T \ \Rightarrow \ \overset{1}{T}=\vp_1^{-1}(T)
\ee
then we may use the inverse function
\bdis
\vp_1^{-1}(T)
\edis
to generate reverse iterations. We have (see (\ref{6.2})) that
\be \label{6.3}
\begin{split}
 & \vp_1(\overset{2}{T})=\overset{1}{T} \ \Rightarrow \ \overset{2}{T}=\vp_1^{-1}(\overset{1}{T})=
 \vp_1^{-1}(\vp_1^{-1}(T))=\vp_1^{-2}(T), \\
 & \vdots \\
 & \overset{k}{T}=\vp_1^{-k}(T),\ k=1,\dots,k_0.
\end{split}
\ee
Of course, we have (see (\ref{6.1}), (\ref{6.3})) that
\bdis
\vp_1^k(\overset{k}{T})=\vp_1^k(\vp_1^{-k}(T))=\vp_1^0(T)=T.
\edis

\subsection{}

Next, the following holds true.

\begin{mydef82}
If
\bdis
t\in [\vp_1^{-k}(T),\vp_1^{-k}(T+H)]
\edis
then (see (\ref{6.3}))
\bdis
\begin{split}
 & \vp_1^r(t)\in [\vp_1^r(\vp_1^{-k}(T)),\vp_1^r(\vp_1^{-k}(T+H))]= \\
 & = [\vp_1^{r-k}(T),\vp_1^{r-k}(T+H)]
\end{split}
\edis
(comp. (\ref{2.10})), i. e.
\be \label{6.4}
\begin{split}
 & \vp_1^0(t)=t\in [\vp_1^{-k}(T),\vp_1^{-k}(T+H)]=[\overset{k}{T},\overset{k}{\wideparen{T+H}}], \\
 & \vp_1^1(t)\in [\vp_1^{-k+1}(T),\vp_1^{-k+1}(T+H)]=[\overset{k-1}{T},\overset{k-1}{\wideparen{T+H}}], \\
 & \vdots \\
 & \vp_1^{k-1}(t)\in [\vp_1^{-1}(T),\vp_1^{-1}(T+H)]=[\overset{1}{T},\overset{1}{\wideparen{T+H}}], \\
 & \vp_1^{k}(t)\in [\vp_1^{0}(T),\vp_1^{0}(T+H)]=[T,T+H].
\end{split}
\ee
\end{mydef82}

\begin{remark}
Of course, the following holds true (see (\ref{5.13}), (\ref{6.4}))
\be \label{6.5}
\begin{split}
 & [T,T+H]\prec [\vp_1^{-1}(T),\vp_1^{-1}(T+H)]\prec \dots \\
 & \prec [\vp_1^{-k}(T),\vp_1^{-1}(T+H)],\ k=1,\dots,k_0.
\end{split}
\ee
Let us remind for comparison that (see \cite{6}, (2.5), Remark 7)
\bdis
\begin{split}
 & [T,T+H]\succ [\vp_1(T),\vp_1(T+H)]\succ \dots \succ \\
 & \succ [\vp_1^k(T),\vp_1^k(T+H)],
\end{split}
\edis
where the direct iteration $\vp_1^k(T)$ is generated by the function $\vp_1(T)$.
\end{remark}

\begin{remark}
The first two reverse iterations of the segments
\bdis
[T,T+\frac{1}{\ln T}],\ [T,T+\frac{1}{\ln^2T}]
\edis
are included in integral formulae (\ref{4.7}).
\end{remark}

\section{Main lemma}

\subsection{}

The following Lemma holds true.

\begin{lemma}
If
\be \label{7.1}
H=o\left(\frac{T}{\ln T}\right),\ T\to\infty
\ee
then for every Lebesgue-integrable function
\bdis
f(t),\ t\in [T,T+H]
\edis
we have that
\be \label{7.2}
\begin{split}
 & \int_T^{T+H} f(t){\rm d}t=\int_{\overset{k}{T}}^{\overset{k}{\wideparen{T+H}}}f[\vp_1^k(t)]
 \prod_{r=0}^{k-1}\tilde{Z}^2[\vp_1^r(t)]{\rm d}t, \\
 & T\to\infty,\ k=1,\dots,k_0;\ \vp_1^0(t)=t,
\end{split}
\ee
where $k_0\in\mbb{N}$ is arbitrary fixed number.
\end{lemma}

\begin{proof}
In our paper \cite{5}, (9.2), (9.5) we have proved the following lemma: if (comp. (\ref{2.3}))
\bdis
\vp_1\{[\overset{1}{T},\overset{1}{\wideparen{T+H}}]\}=[T,T+H]
\edis
then for every Lebesgue-integrable function
\bdis
f(t),\ t\in [T,T+H]
\edis
we have (comp. (\ref{1.2}), (\ref{1.3}))
\be \label{7.3}
\begin{split}
 & \int_T^{T+H} f(t){\rm d}t=\int_{\overset{1}{T}}^{\overset{1}{\wideparen{T+H}}}f[\vp_1(t)]
 \tilde{Z}^2(t){\rm d}t, \\
 & T\geq T_0[\vp_1],\ H\in \left.\left( 0, \frac{T}{\ln T}\right.\right].
\end{split}
\ee
Another form of (\ref{7.3}) is expressed by the formula (see (\ref{1.3}))
\be \label{7.4}
\int_T^{T+H} f(t){\rm d}t=\int_{\overset{1}{T}}^{\overset{1}{\wideparen{T+H}}}f[\vp_1(t)]
\frac{Z^2(t)}{\omega(t)}{\rm d}t.
\ee
Now, the repeated application of the formula (\ref{7.3}) (see (\ref{5.9})) gives the following: if
\bdis
H=o\left(\frac{T}{\ln T}\right)
\edis
then
\bdis
\begin{split}
 & \int_T^{T+H} f(t){\rm d}t=\int_{\overset{1}{T}}^{\overset{1}{\wideparen{T+H}}}f[\vp_1(t)]
 \tilde{Z}^2(t){\rm d}t= \\
 & =\int_{\overset{2}{T}}^{\overset{2}{\wideparen{T+H}}}f[\vp_1^2(t)]
 \tilde{Z}^2[\vp_1(t)]\tilde{Z}^2(t){\rm d}t= \dots = \\
 & = \int_{\overset{k}{T}}^{\overset{k}{\wideparen{T+H}}}f[\vp_1^k(t)]
 \prod_{r=0}^{k-1} \tilde{Z}^2[\vp_1^r(t)]{\rm d}t,
\end{split}
\edis
that is exactly (\ref{7.2}).
\end{proof}

\subsection{}

\begin{remark}
The formula (\ref{7.2}) can be expressed as follows (see (\ref{6.3}))
\bdis
\int_T^{T+H} f(t){\rm d}t=
\int_{\vp_1^{-k}(T)}^{\vp_1^{-k}(T+H)}f[\vp_1^k(t)]\prod_{r=0}^{k-1} \tilde{Z}^2[\vp_1^r(t)]{\rm d}t, \
T\to\infty .
\edis
\end{remark}

\section{Proof of Theorem}

Let
\bdis
\{ f_n(t)\}_{n=1}^\infty,\ t\in [0,2l],\ l=o\left(\frac{T}{\ln T}\right)
\edis
be arbitrary fixed $L_2$-orthogonal system, i. e.
\be \label{8.1}
\int_0^{2l} f_m(t)f_n(t){\rm d}t=
\left\{\begin{array}{rcl} 0 & , & m\not=n, \\ A_n & , & m=n, \end{array}\right. \
A_n=\int_0^{2l} f_n^2(t){\rm d}t .
\ee
Then we have for corresponding system (\ref{2.2}) by (\ref{7.2}), (\ref{8.1}) that
\bdis
\begin{split}
 & f(t) \longrightarrow f_m(\vp_1^k(t)-T)f_n(\vp_1^k(t)-T), \\
 & \int_{\overset{k}{T}}^{\overset{k}{\wideparen{T+2l}}}
 f_m(\vp_1^k(t)-T)f_n(\vp_1^k(t)-T)\prod_{r=0}^{k-1} \tilde{Z}^2[\vp_1^r(t)]{\rm d}t= \\
 & =\int_T^{T+2l}f_m(t-T)f_n(t-T){\rm d}t=  \\
 & =\int_0^{2l}f_m(t)f_n(t){\rm d}t=
 \left\{\begin{array}{rcl} 0 & , & m\not=n, \\ A_n & , & m=n, \end{array}\right.
\end{split}
\edis
i. e. (\ref{2.3}) holds true. Finally, the properties (\ref{2.5}) -- (\ref{2.7}) follows from
(\ref{5.9}), (\ref{5.12}), (\ref{5.13}).

\thanks{I would like to thank Michal Demetrian for his help with electronic version of this paper.}

\end{document}